\documentclass[a4paper,12pt,reqno]{article}
\usepackage{amssymb,amsmath,amsthm}
\usepackage{fullpage}
\usepackage{graphicx}
\newcommand{\iitem}{\vspace{-1.2ex}\item}
\theoremstyle{definition}

\usepackage{color}
\usepackage[normalem]{ulem}

\newcommand\soutD{\bgroup\markoverwith
{\textcolor{DarkGreen}{\rule[.5ex]{2pt}{1pt}}}\ULon}
\newcommand{\Hm}[1]{\leavevmode{\marginpar{\tiny%
$\hbox to 0mm{\hspace*{-0.5mm}$\leftarrow$\hss}%
\vcenter{\vrule depth 0.1mm height 0.1mm width \the\marginparwidth}%
\hbox to
0mm{\hss$\rightarrow$\hspace*{-0.5mm}}$\\\relax\raggedright #1}}}

\title{\textbf{From Symmetry to Monotonicity}}

\author{
 {\sc Bernd Kawohl and David Krej\v ci\v r\'ik} 
 }

\date{\small 12 March 2018}

\begin{document}

\maketitle


\noindent
Symmetry is a fascinating notion in science,
primarily perhaps due to the fact that the symmetry of geometric objects
evokes the aesthetic and elegant perception of every human being.
More fundamentally, the symmetries of spacetime are closely related
to conservation laws in physics. 
In this note we show how one can use a symmetry of a two-dimensional surface
to prove a monotonicity property of a one-dimensional function.

A couple of years ago Jan Ub\o{}e gave 
an unconventional proof of the fact that a certain mapping~$f(x)$ 
which depends on a positive parameter~$a$ 
is increasing in~$x$ if $a>1$ and decreasing if $a<1$,
see~\cite{U}. 
In fact, $f$ is given by
\begin{equation*}
f(x)=\frac{\phi(ax)-\phi(x)}{x \, \big(\Phi(ax)+\Phi(x)\big)}
\end{equation*}
where $\phi(x)=e^{-x^2}$ is the derivative of $\Phi(x)=\int_0^x\phi(u)\, du$. This behaviour of $f$ is used to describe certain effects in a so-called newsvendor model in operations research \cite{UAJ}.
The author of~\cite{U} reckons that 
\emph{Euler would have liked his proof} 
(perhaps because of transferring the problem 
to cumbersome tasks involving infinite sums), 
but also that \emph{a very short proof is lurking around somewhere}.
We believe that our alternative proof 
offers simpler geometric arguments and is more tractable and appealing.

The proof that we 
present makes repeated use of the facts that:

\begin{enumerate}
\iitem[(a)] 
$\phi$ can be extended as an even function on $\mathbb{R}$;
\iitem[(b)]
the derivative satisfies a recurrent formula, namely
$\phi'(x)=-2x\phi(x)$. 
\end{enumerate}

First we observe that the sign of $f'(x)$ is determined by the sign of
\begin{equation}\label{h-initial}
\begin{aligned}
  h(x)
  =& \ 
  \big(ax\phi'(ax)-x\phi'(x)\big)\big(\Phi(ax)+\Phi(x)\big)
  \\
  &\
  -\big(\phi(ax)-\phi(x)\big)\big(\Phi(ax)+\Phi(x)\big)
  \\
  &\ -\big(\phi(ax)-\phi(x)\big)\big(ax\phi(ax)+x\phi(x)\big).
\end{aligned}
\end{equation}
Then we rewrite each of the sums of two terms as integrals
with help of properties~(a) and~(b) above
and arrive at a double-integral representation
\begin{equation}\label{h-final}
  h(x)
  =\int_{-x}^{ax}\int_{-x}^{ax} \phi(u)\phi(v)\  \Gamma(u,v)\ dudv
\end{equation}
with
$$
  \Gamma(u,v) = 2(u+v)(u-v)^2 .
$$

\begin{figure}[h]
\begin{center}
\includegraphics[width=0.5\textwidth]{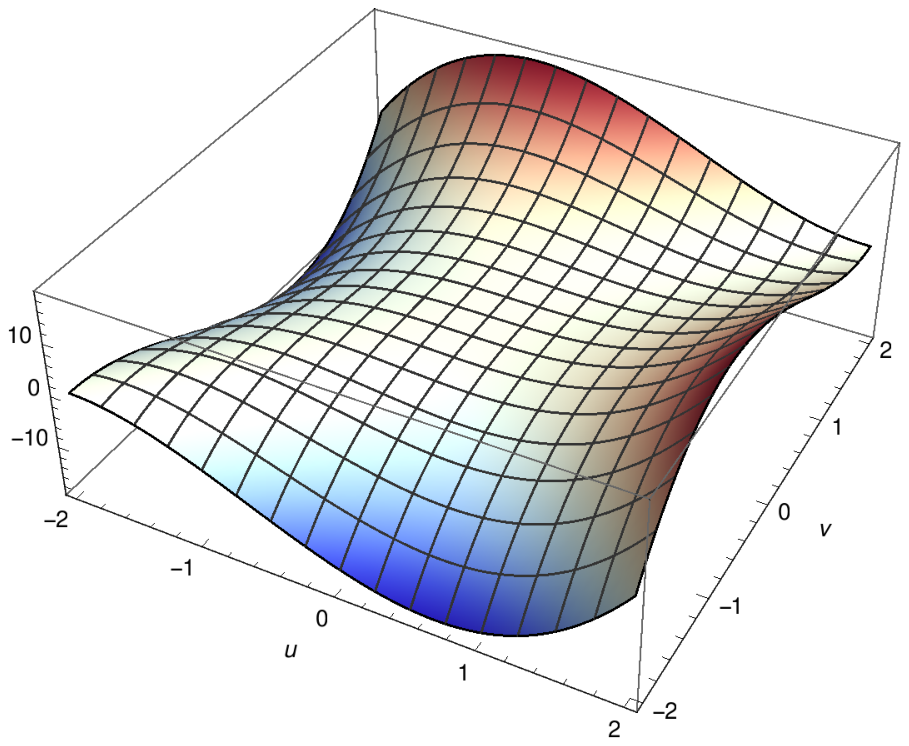}
\
\includegraphics[width=0.38\textwidth]{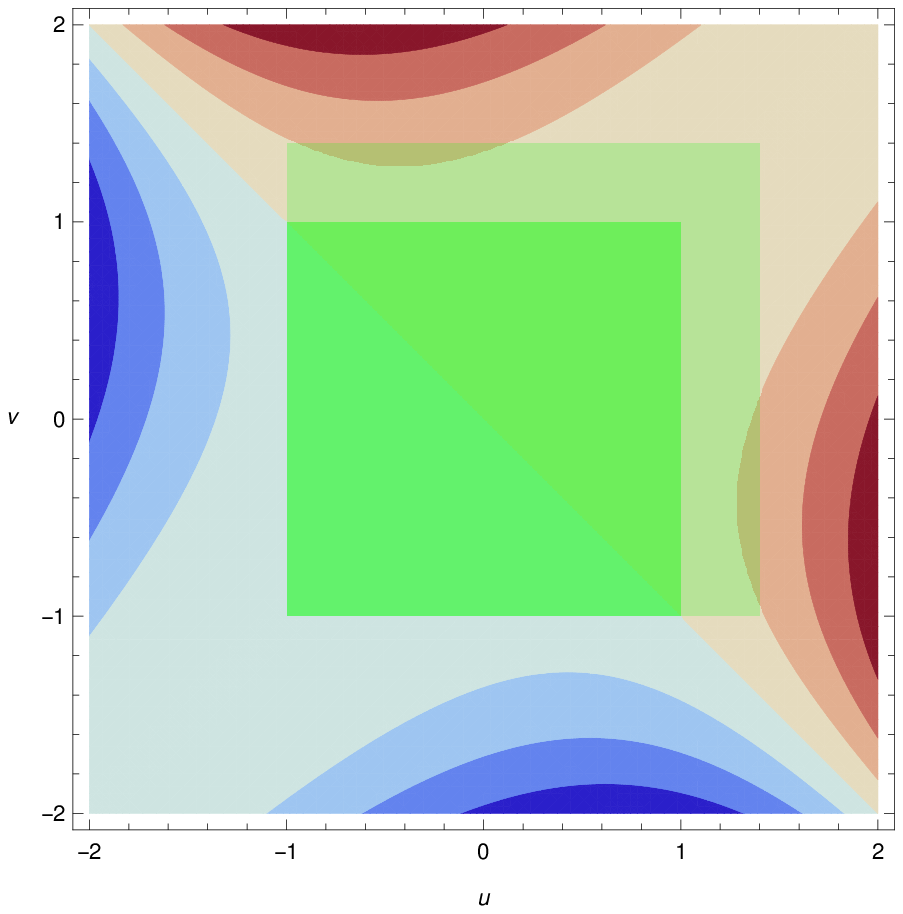}
\includegraphics[width=0.5cm]{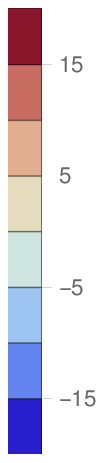}%
\end{center}
\caption{Graph of the surface $\Gamma$ and its contour plot.
Integration over the green square 
(corresponding to $x=1$ and $a=1$ in the picture)
gives zero because of the odd symmetry of~$\Gamma$
with respect to the line $u=-v$,
while the integral over the rectangle expanded 
by the lighter green colour 
into the positive (brownish) region of~$\Gamma$ 
(corresponding to $a=1.4$ in the picture)
gives a positive number.}\label{Fig}
\end{figure}

Notice that the integrand in the last double integral
is odd with respect to the diagonal $v=-u$ and positive above it.
It is due to the obvious symmetry of~$\phi$ stated in~(a) above 
together with a hidden symmetry encoded in the convolution kernel
$$
  \Gamma(u,v)=-\Gamma(-u,-v)
$$
valid for all $u,v \in \mathbb{R}$.
Consequently, since we are integrating over the rectangle
$(-x,ax) \times (-x,ax)$, the odd symmetry immediately implies that
(\emph{cf.}~Figure~\ref{Fig}):
\begin{itemize}
\iitem
for $a=1$ the domains of positivity and negativity of~$\Gamma$
are equilibrated, so $h(x)$ equals zero
(and~$f$ is constant, in fact it is identically equal to zero);
\iitem
for $a>1$ the  domain of positivity of $\Gamma$ outweighs the one of negativity, 
so $h(x)$ is positive
(and~$f$ is strictly increasing); 
\iitem
for $a<1$ the  domain of negativity of $\Gamma$ outweighs the one of positivity, 
so $h(x)$ is negative
(and~$f$ is strictly decreasing). 
\end{itemize}

In summary, we have obtained the desired monotonicity 
of the one-dimensional function~$f$ as a consequence
of the symmetry properties of the  two-dimensional surface~$\Gamma$.

In the rest of this note we explain 
the passage from~\eqref{h-initial} to~\eqref{h-final}.
For the first term appearing in~\eqref{h-initial} we have
\begin{align*}
ax\phi'(ax)-x\phi'(x)
&=-2 \, \big(a^2x^2\phi(ax)-x^2\phi(x)\big)
=-2\int_x^{ax}\big(u^2\phi(u)\big)'\ du\\
&=-4\int_x^{ax}u\phi(u)\ du+4\int_x^{ax}u^3\phi(u)\ du\\
&=-4\int_{-x}^{ax}u\phi(u)\ du+4\int_{-x}^{ax}u^3\phi(u)\ du.
\end{align*}
In the last equality we have used the fact that $u\phi(u)$ and $u^3\phi(u)$ are odd functions of $u$.
Similarly we write
\begin{align*}
\Phi(ax)+\Phi(x) &= \int_{-x}^{ax}\phi(u)\ du,
\\
\phi(ax)-\phi(x) &=\int_{-x}^{ax}\phi'(u)\ du=-2\int_{-x}^{ax}u\phi(u)\ du,
\\
ax\phi(ax)+x\phi(x)&=\int_{-x}^{ax}(u\phi(u))'\ du
=\int_{-x}^{ax}\phi(u)\ du-2\int_{-x}^{ax}u^2\phi(u)\ du.
\end{align*}
Putting these formulae together
one can rewrite~$h$ as follows:
\begin{equation*}
h(x)
= 4\left(\int_{-x}^{ax} u^3\phi(u)\ du\right) 
  \left(\int_{-x}^{ax}\phi(u)\ du\right)
-4\left(\int_{-x}^{ax} u\phi(u)\ du\right) 
  \left(\int_{-x}^{ax} u^2\phi(u)\ du\right).
\end{equation*}

Now the main idea is to rewrite the product of integrals 
$(\int g(u)\, du)(\int j(v)\, dv)$ 
as a double integral $\iint g(u)j(v)\, dudv$ and realise that
\begin{equation*}
  h(x)
  =\int_{-x}^{ax}\int_{-x}^{ax} \phi(u)\phi(v)\  \tilde\Gamma(u,v)\ dudv
  =\int_{-x}^{ax}\int_{-x}^{ax} \phi(u)\phi(v)\  \tilde\Gamma(v,u)\ dudv,
\end{equation*}
where $\tilde\Gamma(u,v)=4(u^3-uv^2)$.
Here, the second equality is due to the fact 
that the role of variables~$u$ and~$v$ can be interchanged.
Therefore, $\tilde\Gamma(u,v)$ can be 
replaced by  $\tilde\Gamma(v,u)$ and also
by the arithmetic mean
$$
  \Gamma(u,v) = \frac{\tilde\Gamma(u,v) + \tilde\Gamma(v,u)}{2}
  =2(u+v)(u-v)^2, 
$$
so that~\eqref{h-final} holds.

\subsubsection*{Acknowledgment}
D.K. was partially supported 
by the GACR grant No.\ 18-08835S
and by FCT (Portugal)
through project PTDC/MAT-CAL/4334/2014.


\renewcommand\refname


\subsubsection*{Authors' addresses}

Bernd Kawohl,
Mathematisches Institut,
Universit\"at zu K\"oln,
D-50923 K\"oln,
Germany;
kawohl@mi.uni-koeln.de

\medskip\noindent
David Krej\v ci\v r\' ik, 
Department of Mathematics, 
Faculty of Nuclear Sciences and Physical Engineering, 
Czech Technical University in Prague, 
Trojanova 13, 
120\,00 Prague, 
Czech Republic;
david.krejcirik@fjfi.cvut.cz

\end{document}